\documentclass[12pt]{article}
\usepackage{amsfonts,amssymb,amsmath,amscd,amsthm,graphicx}

\newtheorem{thm}{Theorem}[section]
\newtheorem{lem}[thm]{Lemma}
\newtheorem{prop}[thm]{Proposition}

\newtheorem{prob}[thm]{Problem}
\theoremstyle{definition}
 
\newtheorem{rmk}[thm]{Remark}
\newtheorem{ex}[thm]{Example}

\numberwithin{thm}{section}
\numberwithin{equation}{section}

\newcommand{\Hom}{\operatorname{Hom}}
\newcommand{\Ext}{\operatorname{Ext}}
\newcommand{\Spec}{\operatorname{Spec}}

\newcommand{\ob}{\operatorname{ob}}
\newcommand{\Hilb}{\operatorname{Hilb}}
\newcommand{\red}{\operatorname{red}}
\newcommand{\im}{\operatorname{im}}

\newcommand{\Pic}{\operatorname{Pic}}
\newcommand{\coker}{\operatorname{coker}}

\newcommand{\car}{\operatorname{char}}
\newcommand{\ov}[1]{\overline{#1}}

\renewcommand{\labelenumi}{{\rm (\arabic{enumi})}}

\newcommand{\Proof}{\noindent {\it Proof}.\ \ }
\newcommand{\QED}{$\Box$ \vskip 3mm}

\newcommand{\mapright}[1]{%
\smash{\mathop{%
\hbox to 1cm{\rightarrowfill}}\limits^{#1}}}
\newcommand{\mapleft}[1]{%
\smash{\mathop{%
\hbox to 1cm{\leftarrowfill}}\limits^{#1}}}
\newcommand{\mapdown}[1]{\Big\downarrow
\rlap{$\vcenter{\hbox{$\scriptstyle#1\,$}}$ }}
\newcommand{\mapup}[1]{\Big\uparrow
\rlap{$\vcenter{\hbox{$\scriptstyle#1$}}$ }}
\newcommand{\mapid}{\big\Vert
\rlap{$\vcenter{\hbox{$\,\scriptstyle {\rm id}$}}$ }}

\newcommand{\smapup}[1]{\uparrow
\rlap{$\vcenter{\hbox{$\scriptstyle#1$}}$ }}

\newcommand{\rmapdown}[1]{\stackrel{|}{\downarrow}
\rlap{$\vcenter{\hbox{$\scriptstyle#1\,$}}$ }}

\title{Obstructions to deforming curves on a $3$-fold, I: \\
{\Large A generalization of Mumford's example\\
and an application to Hom schemes}}
\author{Shigeru {\sc Mukai}\thanks{Supported in part by 
the JSPS Grant-in-Aid for Scientific Research (B) 17340006.}\,
and Hirokazu {\sc Nasu}\thanks{Supported in part by the 
21st century COE program
``Formation of an International Center of Excellence
in the Frontier of Mathematics
and Fostering of Researchers in Future Generations''.}
}
\date{}

\begin{document}
\maketitle

\begin{abstract}
We give a sufficient condition for a first order infinitesimal 
deformation of a curve on a 3-fold to be obstructed.
As application we construct generically non-reduced components 
of the Hilbert schemes of uniruled 3-folds and the Hom scheme 
from a general curve of genus five to a general cubic 3-fold.
\end{abstract}

\section{Introduction}

We study the (embedded) deformation of a (smooth projective) curve 
$C$ on a smooth projective $3$-fold $V$ under the presence of 
a certain pair of a smooth surface $S$ and a smooth curve $E$ 
such that $C, E \subset S \subset V$.
In other words we study the Hilbert scheme $\Hilb^{sc} V$ 
of smooth curves on $V$ with the help of intermediate surfaces.
Let $\tilde C$ be a first order infinitesimal deformation of $C \subset V$.
As is well known, $\tilde C$ determines a global section $\alpha$ 
of the normal bundle $N_{C/V}$.
It also determines a cohomology class $\ob(\alpha) \in H^1(N_{C/V})$ 
such that $\tilde C$ lifts to a deformation over $\Spec k[t]/(t^3)$
if and only if $ \ob(\alpha)$
is zero (\S\ref{general theory}). 
This $\ob(\alpha)$ is called the {\em (primary) obstruction} 
of $\alpha$ (or $\tilde C$).
It is generally difficult to compute $\ob(\alpha)$ for given $\alpha$.
In this paper, we give a sufficient condition for $\ob(\alpha)\ne 0$ 
in terms of $\pi_S(\alpha)$, the {\em exterior component} of $\alpha$
(Theorem~\ref{thm:main thm1}).
Here $\pi_S$ is the natural projection 
$N_{C/V} \rightarrow N_{S/V}\big{\vert}_C$.

In each of the following examples, 
the tangent space $t_{W,C}$ of the subvariety $W$ is everywhere 
of codimension one in $ H^0(N_{C/V})$, the tangent space of $\Hilb^{sc} V$.
If $[C]\in W$ is general, then every $\alpha \in H^0(N_{C/V})$ 
not in $t_{W,C}$ satisfies the condition, hence is obstructed.
Therefore $\Hilb^{sc} V$ is everywhere non-reduced along $W$.

\begin{ex}\label{ex:pathology}
 Let $V$ be the projective space $\mathbb P^3$, 
 $S \subset \mathbb P^3$ a smooth cubic surface, $E \subset S$ 
 a $(-1)$-$\mathbb P^1$ and $C \subset S$ a smooth member
 of the linear system $|4h+2E|\simeq \mathbb P^{37}$ on $S$. 
 $C$ is of degree $14$ and genus $24$.
 Such $C$'s are parametrized by 
$W = W^{56} \subset \Hilb^{sc} \mathbb P^3$, an open subset
 of a $\mathbb P^{37}$-bundle over $|3H|\simeq \mathbb P^{19}$.
 Here $H$ is a plane in $\mathbb P^3$ and $h$ is its restriction to $S$.
\end{ex}

\begin{ex}\label{ex:cubic 3-fold}
 Let $V$ be a smooth cubic $3$-fold $V_3 \subset \mathbb P^4$,
 $S$ its general hyperplane section, $E \subset S$ 
 a $(-1)$-$\mathbb P^1$ and $C \subset S$ a smooth member of
 $|2h+2E|\simeq \mathbb P^{12}$. 
 $C$ is of degree $8$ and genus $5$.
Such $C$'s are parametrized by $W=W^{16} \subset \Hilb^{sc} V$, 
an open subset of $\mathbb P^{12}$-bundle over the dual projective 
space $\mathbb P^{4,\vee}$.
Here $h$ is the restriction to $S$ of a hyperplane $H$ of $\mathbb P^4$.
(See \S\ref{dichotomy} for details.)
\end{ex}

For many uniruled 3-folds $V$, we can find a curve $C \subset V$ 
similar to the above examples.
More precisely we have the following:
 
\begin{thm}\label{thm:main thm2}
 Suppose that $E$ is a $(-1)$-$\mathbb P^1$ on $S$, $N_{E/V}$ 
 is generated by global sections and $p_g(S)=h^1(N_{S/V})=0$.
 Then the Hilbert scheme $\Hilb^{sc} V$ of smooth curves on $V$
 contains infinitely many generically non-reduced 
 (irreducible) components 
 $W_n$ $(n\in \mathbb Z_{\ge 0})$ with the following property:
 \begin{enumerate}
 \renewcommand{\labelenumi}{{\rm (\alph{enumi})}}
  \item every member of $W_n$ is contained in a deformation of $S$ in $V$, 
	and
  \item every general member of $W_n$ has a first order 
	infinitesimal deformation whose primary obstruction is nonzero.
 \end{enumerate}
\end{thm}

See Example~\ref{ex:uniruled} for 3-folds $V$ with such $S$ and $E$.

Mumford~\cite{Mumford} proved the non-reducedness of 
$\Hilb^{sc} \mathbb P^3$ (Example~\ref{ex:pathology}) 
by a global argument but later 
Curtin~\cite{Curtin} gave another proof by 
infinitesimal analysis of deformations.
Recently Nasu~\cite{Nasu} has simplified and generalized Curtin's proof.
This theorem follows the line of these works.
Vakil~\cite{V} has also shown that various moduli spaces have
non-reduced components by a different method
({\it cf}. Remark~\ref{rmk:murphy}).

For a given projective scheme $X$, 
the set of morphisms $f:X\rightarrow V$ has a natural scheme 
structure as a subscheme of the Hilbert scheme of $X \times V$. 
This scheme is called the {\em Hom scheme} and denoted by $\Hom(X,V)$.
When we fix a projective embedding $V \hookrightarrow \mathbb P^n$,
all the morphisms of degree $d$ are parametrized by an open and 
closed subscheme, which we denote by $\Hom_d(X,V)$. 
Our Example~\ref{ex:cubic 3-fold} gives rise to a counterexample 
to the following problem on the Hom scheme:

\begin{prob}[$k=\mathbb C$]\label{hom version}
Is every component of $\Hom(X,V')$ generically smooth for a smooth 
curve $X$ with general modulus and a general member $V'$ 
of the Kuranishi family of $V$?
\end{prob}

Let $\Hom_8(X,V_3)$ be the Hom scheme of morphisms of degree $8$ 
from a general curve $X$ of genus $5$ to a smooth cubic $3$-fold 
$V_3 \subset \mathbb P^4$.
Its expected dimension equals 4 ({\it cf}. \S\ref{hom scheme}).

\begin{thm}[$\car k=0$]\label{thm:main thm}
Assume that $V_3$ is either (moduli-theoretically) general or of
 Fermat type $$V_3^{\rm Fermat}: x_0^3+x_1^3+x_2^3+x_3^3+x_4^3=0.$$ 
Then $\Hom_8(X,V_3)$ has an irreducible component 
of expected dimension which is generically non-reduced.
\end{thm}

In order to prove Theorem~\ref{thm:main thm2}, 
we take a rational section $v$ of the normal bundle $N_{S/V}$. 
Suppose that $v$ has a pole only along a smooth curve $E \ne C$ 
and of order one, that is, $v \in H^0(N_{S/V}(E)) \setminus H^0(N_{S/V})$.
Then the divisor $(v)_0$ of zero does not contain $E$ as a component.
The restriction $v\big{\vert}_C$ belongs to 
$H^0(N_{S/V}\big{\vert}_C) \subset H^0(N_{S/V}(E)\big{\vert}_C)$ 
if and only if
\begin{equation}\label{ineq:intersection}
(v)_0 \cap E \ge C \cap E 
\end{equation}
as a divisor on $E$. 

\begin{thm}\label{thm:main thm1}
 Let $C,E \subset S \subset V$ and 
 $v \in H^0(N_{S/V}(E)) \setminus H^0(N_{S/V})$ be as above and 
 assume that $(E^2) < 0$.
 If the following conditions are satisfied, then
  every first order infinitesimal deformation $\tilde C$ of 
 $C \subset V$, or $\alpha$, whose exterior
 component coincides with $v\big{\vert}_C$ is obstructed.
 \begin{enumerate}
 \renewcommand{\labelenumi}{{\rm (\alph{enumi})}}
 \item The equality holds in \eqref{ineq:intersection}.
 \item Let $\partial$ denote the coboundary map of the 
	exact sequence of
	\begin{equation}\label{ses:normal bundle of line}
	 0 \longrightarrow N_{E/S}
	 \longrightarrow N_{E/V}
	 \longrightarrow N_{S/V}\big{\vert}_E
	 \longrightarrow 0 
	\end{equation}
	tensored with $\mathcal O_S(E)$.
	Then the image $\partial(v\big{\vert}_E)$ of 
	$v\big{\vert}_E \in H^0(N_{S/V}(E)\big{\vert}_E)$ 
	is nonzero in $H^1(N_{E/S}(E))\simeq H^1(\mathcal O_E(2E))$.
 \item The restriction map 
	$H^0(S,\Delta) \rightarrow H^0(E,\Delta\big{\vert}_E)$
	is surjective, where $\Delta:=C-2E+K_V\big{\vert}_S$.
 \end{enumerate}
\end{thm}

If $E \subset S$ is a $(-1)$-$\mathbb P^1$, $N_{E/V}$ is generated by
global sections and $v\big{\vert}_E$ is a general member of 
$H^0(N_{S/V}(E)\big{\vert}_E)$, then the condition (b) 
is satisfied (Lemma~\ref{lem:nonzero}).
If $E$ is a $(-1)$-$\mathbb P^1$, we have $(\Delta. E)=0$ and hence
(c) is equivalent to that $|\Delta -E|+E \ne |\Delta|$.

We prove Theorem~\ref{thm:main thm1} in \S\ref{obstruction} 
and Theorem~\ref{thm:main thm2} in \S\ref{hilbertscheme}. 
In the final section, we prove Theorem~\ref{thm:main thm}.

We work over an algebraically closed field $k$ in arbitrary
characteristic except for in \S\ref{hom scheme}.
We denote by $(A. B)$ the intersection number of two divisors $A$ 
and $B$ on a surface.
For a subscheme $S$ of $V$ and a sheaf $\mathcal F$ on $V$, we denote 
the {\em restriction map} 
$H^i(V, \mathcal F) \to H^i(S, \mathcal F\big{\vert}_S)$ 
by $\big{\vert}_S$.

\section{Obstruction to deforming curves}\label{obstruction}

\subsection{General theory}\label{general theory}

Let $C$ be a smooth closed subvariety of a smooth variety $V$.
We denote the normal bundle $\mathcal Hom(\mathcal I_C,\mathcal O_C)$ 
of $C$ in $V$ by $N_{C/V}$,
where $\mathcal I_C$ is the ideal sheaf of $C$ in $V$.
An {\em (embedded) first order infinitesimal deformation} of
$C \subset V$ 
is a closed subscheme 
$\tilde C \subset V\times \Spec k[t]/(t^2)$ 
which is flat over $\Spec k[t]/(t^2)$ and whose central fiber is 
$C \subset V$.
Let $\mathcal I_{\tilde C}$ be the ideal sheaf of $\tilde C$, 
which is also flat over $k[t]/(t^2)$.
The multiplication endomorphism of $\mathcal O_V \otimes k[t]/(t^2)$ 
by $t$ induces a homomorphism 
$\mathcal I_{\tilde C} \to \mathcal O_{\tilde C}$, 
which factors through 
$\alpha: \mathcal I_C \to t\mathcal O_C \simeq \mathcal O_C$.
Moreover, $\tilde C$ is recovered from the homomorphism 
$\alpha: \mathcal I_C \to \mathcal O_C$.
In the sequel we identify $\tilde C$ with $\alpha \in H^0(N_{C/V})$.

The standard exact sequence
\begin{equation}\label{ses:standard}
0 \longrightarrow \mathcal I_C
 \longrightarrow \mathcal O_V
 \longrightarrow \mathcal O_C
 \longrightarrow 0
\end{equation}
induces 
$\delta: H^0(N_{C/V}) =\Hom(\mathcal I_C,\mathcal O_C)
\rightarrow \Ext^1(\mathcal I_C,\mathcal I_C)$
as a coboundary map.
Then $\tilde C$ lifts to a deformation over $\Spec k[t]/(t^3)$
if and only if 
 $$
 \ob(\alpha):= \delta(\alpha) \cup \alpha 
 = \alpha \cup \mathbf k_{C,V} \cup \alpha
 \in \Ext^1(\mathcal I_C,\mathcal O_C)
 $$
is zero, where 
$\mathbf k_{C,V} \in \Ext^1(\mathcal O_C,\mathcal I_C)$ 
is the extension class of \eqref{ses:standard}. 
$\ob(\alpha)$ is called the {\em obstruction} of $\alpha$ (or
$\tilde C$). Since both $C$ and $V$ are smooth, $\ob(\alpha)$ is contained
in the subspace $H^1(N_{C/V})\subset \Ext^1(\mathcal I_C,\mathcal O_C)$
({\it cf}.~\cite[Chap.~I, Proposition 2.14]{Kollar}).

The tangent space of the Hilbert scheme $\Hilb V$ at $[C]$
is isomorphic to $H^0(N_{C/V})$.
If $\Hilb V$ is nonsingular at $[C]$,
then every first order infinitesimal deformation of $C \subset V$ 
lifts to a deformation over $\Spec k[t]/(t^n)$ for any $n \ge 3$.
If $\ob(\alpha)\ne 0$ for some $\alpha\in H^0(N_{C/V})$,
then $\Hilb V$ is singular at $[C]$.

Let $L$ be a line bundle on $V$ and 
$
\delta: H^0(C, L\big\vert_C)\rightarrow 
H^1(V, L \otimes \mathcal I_C)
$, $u \mapsto u \cup \mathbf k_{C,V}$, 
the coboundary map of the exact sequence $L \otimes$\eqref{ses:standard}.
We denote the composite of $\delta$ 
and the restriction map 
$
\big{\vert}_C: 
H^1(V,L \otimes \mathcal I_C) \rightarrow 
H^1(C,L\big\vert_C \otimes {N_{C/V}}^{\vee})
$
by 
\begin{equation}\label{eq:d-map}
d_{C, L}: H^0(C, L\big\vert_C) \longrightarrow 
H^1(C,L\big\vert_C \otimes {N_{C/V}}^{\vee}).
\end{equation}
If $V$ is projective, $C$ is a divisor and $L=\mathcal O_V(C)$, then
\begin{equation}\label{eq:divisorial d-map}
d_{C, \mathcal O_V(C)}: H^0(C, N_{C/V}) \longrightarrow H^1(C,\mathcal O_C)
\end{equation}
is the tangential map of the natural morphism $C' \mapsto \mathcal
O_C(C')$ from the Hilbert scheme of divisors $C' \subset V$
to the Picard scheme $\Pic C$.

\subsection{Exterior component}\label{exterior} 

From now on we assume that $C$ is contained in a smooth divisor $S \subset V$.
There exists a natural exact sequence
\begin{equation}\label{ses:normal bundle of curve}
 0 \longrightarrow N_{C/S}
 \longrightarrow N_{C/V}
 \stackrel{\pi_S}\longrightarrow N_{S/V}\big{\vert}_C
 \longrightarrow 0
\end{equation}
of normal bundles.
In this article we compute not $\ob(\alpha)$ itself but its image by 
$$H^1(\pi_S): H^1(N_{C/V}) \longrightarrow H^1(N_{S/V}\big{\vert}_C).$$
We call the image 
the {\em exterior component} of $\ob(\alpha)$ and denote by $\ob_S(\alpha)$.
Here we give another expression of $\ob_S(\alpha)$.
Let 
$$
d_{C, \mathcal O_V(S)}:
H^0(N_{S/V}\big\vert_C) \simeq H^0(\mathcal O_C(S)) 
\longrightarrow H^1(C,{N_{C/V}}^{\vee}\otimes N_{S/V}\big{\vert}_C)
$$ 
be the map \eqref{eq:d-map} for the line bundle $L= \mathcal O_V(S)$.
We abbreviate this as $d_C$.

\begin{lem}\label{lem:exterior component}
 \begin{equation*}
 \ob_S(\alpha)= d_C(\pi_S(\alpha)) \cup \alpha,
 \end{equation*}
where $\cup$ is the cup product map
 \begin{equation}\label{map:cup exterior}
 H^1(C,{N_{C/V}}^{\vee}\otimes N_{S/V}\big{\vert}_C) \times H^0(C,N_{C/V})
 \overset{\cup}{\longrightarrow} H^1(N_{S/V}\big{\vert}_C).
 \end{equation}
 \end{lem}
 \Proof
For each $i=0$ and $1$, $H^i(\pi_S)$ is equal to the restriction
to $H^i(N_{C/V}) \subset \Ext^i(\mathcal I_C,\mathcal O_C)$ of
the cup product map 
$\Ext^i(\mathcal I_C,\mathcal O_C)
\overset{\cup\, \iota}{\longrightarrow} \Ext^i(\mathcal I_S,\mathcal O_C)$,
where $\iota: \mathcal I_S \hookrightarrow \mathcal I_C$ 
is the natural inclusion.
Recall that the coboundary map $\delta$ in \S\ref{general theory}
is also a cup product map with the extension class
$\mathbf k_{C,V}$ of \eqref{ses:standard}.
Therefore we have
 $$
 \pi_S(\ob(\alpha))
 =\iota \cup (\alpha \cup \mathbf k_{C,V} \cup \alpha)
 = \pi_S(\alpha) \cup \mathbf k_{C,V} \cup \alpha
 = d_C (\pi_S(\alpha)) \cup \alpha. 
 \qquad \Box
$$

Let $d_S: H^0(C, N_{S/V}) \rightarrow H^1(S,\mathcal O_S)$ 
be the map \eqref{eq:divisorial d-map} for $S \subset V$.
$d_S$ and $d_C$ are closely related by the
following commutative diagram:
\begin{equation}\label{diag:d-map}
\begin{array}{ccccc}
H^0(N_{S/V}) & \mapright{d_S} & 
H^1(\mathcal O_S) \\
 & & \mapdown{|_C}\\
 \mapdown{|_C} && H^1(\mathcal O_C) \\
 && \mapdown{H^1(\iota)} \\
H^0(N_{S/V}\big{\vert}_C) & \mapright{d_C} & 
 H^1({N_{C/V}}^{\vee}\otimes N_{S/V}),
\end{array}
\end{equation}
 where $\iota: \mathcal O_C \rightarrow {N_{C/V}}^{\vee}\otimes N_{S/V}$ 
 is the inclusion induced by $\pi_S$ of \eqref{ses:normal bundle of curve}.
 In some cases the exterior component of $\ob(\alpha)$ depends only on that of $\alpha$.

\begin{lem}\label{lem:exterior component2}
 If $\pi_S(\alpha) \in H^0(N_{S/V}\big{\vert}_C)$ 
 lifts to a global section $v$ of $N_{S/V}$, then we have 
 \begin{equation*}
 \ob_S(\alpha)= d_S(v)\big{\vert}_C \cup \pi_S(\alpha).
 \end{equation*}
 \end{lem}
\Proof
By the diagram \eqref{diag:d-map} we have
$d_C (\pi_S(\alpha))=H^1(\iota)(d_S(v)\big{\vert}_C)$.
By the commutative diagram
$$
\begin{array}{ccccc}
d_C(\pi_S(\alpha)) && \alpha && \\ [-10pt]
\rotatebox{-90}{$\in$} && \rotatebox{-90}{$\in$} && \\ [4pt]
H^1({N_{C/V}}^{\vee}\otimes N_{S/V}) & \times 
 & H^0(N_{C/V}) & \overset{\cup}{\longrightarrow} 
& H^1(N_{S/V}\big{\vert}_C)\\
\mapup{H^1(\iota)} && \mapdown{\pi_S} && \Vert \\
H^1(\mathcal O_C) & \times & H^0(N_{S/V}\big{\vert}_C) 
 & \overset{\cup}{\longrightarrow} & H^1(N_{S/V}\big{\vert}_C), \\
\rotatebox{90}{$\in$} && && \\ [-4pt] 
d_S(v)\big{\vert}_C && && \\ 
\end{array}
\begin{array}{ccc}
\end{array}
$$
whose first cup product map is \eqref{map:cup exterior},
we have 
$$
\ob_S(\alpha) = d_C (\pi_S(\alpha)) \cup \alpha
=d_S(v)\big{\vert}_C \cup \pi_S(\alpha)
$$
in $H^1(N_{S/V}\big{\vert}_C)$ by Lemma~\ref{lem:exterior component}. 
\QED

\subsection{Infinitesimal deformation with a pole}\label{with pole}

We assume that $V$ is a 3-fold, $S \subset V$ is a smooth surface and
$E$ is a smooth curve on $S$ with $(E^2)< 0$ as in Theorem~\ref{thm:main thm1}.
We denote the complemental open varieties 
$S \setminus E$ and $V \setminus E$
by $S^{\circ}$ and $V^{\circ}$, respectively, 
and the map \eqref{eq:divisorial d-map} for $S^{\circ} \subset V^{\circ}$ by
$d_{S^{\circ}}: H^0(N_{S^{\circ}/V^{\circ}}) 
\rightarrow H^1(\mathcal O_{S^{\circ}})$.
In this subsection we study the singularity of
$d_{S^{\circ}}(v) \in H^1(\mathcal O_{S^{\circ}})$ along the 
boundary $E$ for $v \in H^0(S,N_{S/V}(E))$
(an infinitesimal deformation with a pole).
The pole of $d_{S^{\circ}}(v)$ is of order at most $2$ and its principal
part coincides with $\partial(v\big{\vert}_E)$
(Proposition~\ref{prop:key proposition}).

Let $\iota: S^{\circ} \hookrightarrow S$ be the open immersion.
Then $\iota_*\mathcal O_{S^{\circ}}$ contains $\mathcal O_S(nE)$
as a subsheaf for any $n \ge 0$.
There exists a natural inclusion
$
\mathcal O_S \subset \mathcal O_S(E) \subset \cdots
\subset \mathcal O_S(nE) \subset \cdots
$
and $\iota_*\mathcal O_{S^{\circ}}$
is the inductive limit $\lim_{n\rightarrow \infty} \mathcal O_S(nE)$.

\begin{lem}\label{lem:injectivity}
 Let $L$ be a line bundle on $S$.
 If $\deg(L\big{\vert}_E) \le 0$, then
 $$
 H^1(S,L) \longrightarrow H^1(S^{\circ},L\big{\vert}_{S^{\circ}})
 $$
 induced by the inclusion 
 $L \hookrightarrow L\otimes \iota_*\mathcal O_{S^{\circ}}$ is injective.
\end{lem}
\Proof
There exists an open affine finite covering 
$\mathfrak U=\left\{ U_i \right\}_{i=1,\ldots,n}$ of $S$. 
Let $\mathbf c = \{c_{ij}\}_{1 \le i < j \le n}$ be a 
$1$-cocycle with coefficient $L$ with respect to $\mathfrak U$ 
and $\gamma_m$ its cohomology class in $H^1(S,L(mE))$ for every $m \ge 0$.
If $\mathbf c$ is a $1$-coboundary of $L\big{\vert}_{S^{\circ}}$, then 
$\mathbf c$ becomes that of $L(mE)$, that is, $\gamma_m=0$, 
for a sufficiently large $m$. Since $\deg(L\big{\vert}_E)\le 0$, 
we have $H^0(E,L(mE)\big{\vert}_E)=0$ for $m \ge 1$.
Hence
$$
H^1(S,L((m-1)E)) \longrightarrow H^1(S,L(mE)) 
$$
is injective. 
Therefore, $\gamma_{m-1}$ is also $0$.
By induction $\gamma_0$ is zero in $H^1(S,L)$.
\QED

By the lemma, the natural map
$H^1(\mathcal O_S(2E))\rightarrow H^1(\mathcal O_{S^{\circ}})$ is
injective.
We identify $H^0(N_{S/V}(E))$ and $H^1(\mathcal O_S(2E))$
with their images in $H^0(N_{S^{\circ}/V^{\circ}})$ 
and $H^1(\mathcal O_{S^{\circ}})$, respectively.

\begin{prop}\label{prop:key proposition}
\begin{enumerate}
 \item $d_{S^{\circ}}(H^0(S,N_{S/V}(E)))
  \subset H^1(S,\mathcal O_S(2E)).$
 \item Let $d_S$ be the restriction of $d_{S^{\circ}}$
	to $H^0(S,N_{S/V}(E))$ and let
	$\partial$ be the coboundary map in Theorem~\ref{thm:main thm1}.
	Then the diagram
	\begin{equation*}
	 \begin{array}{ccc}
	 H^0(S,N_{S/V}(E)) & \mapright{d_{S}}
	 & H^1(S,\mathcal O_S(2E)) \\ 
	 \mapdown{|_E} && \mapdown{|_E} \\
	 H^0(E,N_{S/V}(E)\big{\vert}_E) & \mapright{\partial} 
	 & H^1(E,\mathcal O_E(2E))
	 \end{array}
	\end{equation*}
	is commutative.
\end{enumerate}
\end{prop}
\Proof
Let $\mathfrak U:=\{ U_i \}_{i \in I}$ be an affine open covering of $V$ 
and let $x_i=y_i=0$ be the local equation of $E$ over $U_i$ such that
$y_i$ defines $S$ in $U_i$.
Through the proof, for a local section $t$ of a sheaf $\mathcal F$ on $V$, 
$\bar t$ denotes the restriction 
$t\big{\vert}_S \in \mathcal F\big{\vert}_S$ for conventions.
Let $D_{x_i}$ and $D_{\bar x_i}$ denote the affine open subsets of 
$U_i$ and $U_i \cap S$ defined by $x_i \ne 0$ and $\bar x_i \ne 0$,
respectively.
Then 
$\left\{ D_{\bar x_i} \right\}_{i \in I}$ is an affine open covering of 
$S^{\circ}$ since $D_{\bar x_i}=D_{x_i} \cap S=U_i \cap S^{\circ}$.

Let $v$ be a global section of $N_{S/V}(E) \simeq \mathcal O_S(S)(E)$.
Then the product $\bar x_i v$ is contained in $H^0(U_i,\mathcal O_S(S))$
and lifts to a section
$s_i' \in \Gamma(U_i, \mathcal O_V(S))$ since $U_i$ is affine.
In particular,
$v$ lifts to the section $s_i:=s_i'/x_i$ of 
$\mathcal O_{V^{\circ}}(S^{\circ})$ over $D_{x_i}$.
Then we have
$$
\delta(v)_{ij}= s_j -s_i
\qquad
\mbox{in}
\qquad
\Gamma(D_{x_i} \cap D_{x_j}, \mathcal O_{V^{\circ}}(S^{\circ}))
$$
for every $i,j$,
where $\delta: H^0(\mathcal O_{S^{\circ}}(S^{\circ}))
\rightarrow H^1(\mathcal O_{V^{\circ}})$ is the coboundary map of
$$
[0 \longrightarrow \mathcal I_{S^{\circ}}
\longrightarrow \mathcal O_{V^{\circ}}
\longrightarrow \mathcal O_{S^{\circ}}
\longrightarrow 0
] \otimes \mathcal O_{V^{\circ}}(S^{\circ})
$$
in \S\ref{general theory}.
Since $v$ is a global section of $\mathcal O_{S^{\circ}}(S^{\circ})$,
$\delta(v)_{ij}$ is contained in 
$\Gamma(D_{x_i} \cap D_{x_j}, \mathcal O_{V^{\circ}})$.
Moreover since $x_is_i=s_i' \in \Gamma(U_i,\mathcal O_V(S))$ 
for every $i$, $f_{ij}:=x_ix_j\delta(v)_{ij}$ is contained in 
$\Gamma(U_i \cap U_j, \mathcal O_V)$.
Hence we have
$$
d_{S^{\circ}}(v)_{ij} =(\delta(v)_{ij})\big{\vert}_{S^\circ}= 
\dfrac{\bar f_{ij}}{\bar x_i \bar x_j}
\qquad
\mbox{in}
\qquad
\Gamma(D_{\bar x_i} \cap D_{\bar x_j}, \mathcal O_{S^{\circ}}), 
$$
where $\bar f_{ij}$ is the restriction of 
$f_{ij} \in \mathcal O_{U_i \cap U_j}$ to $S \cap U_i \cap U_j$.
By definition $\bar f_{ij}$ 
belongs to $\Gamma(S\cap U_i \cap U_j, \mathcal O_S)$.
Hence $d_{S^{\circ}}(v)_{ij}$ is contained in 
$\Gamma(S\cap U_i \cap U_j, \mathcal O_S(2E))$.
Thus we have proved (1).

Now we compute the image of $d_S(v)=d_{S^{\circ}}(v)$ by
the restriction map 
$H^1(S,\mathcal O_S(2E)) \rightarrow H^1(E,\mathcal O_E(2E))$
regarding $\mathcal O_E(2E)$ 
as the quotient sheaf $\mathcal O_S(2E)/\mathcal O_S(E)$.
For the computation,
we need to consider the relation between the local equations 
$x_i=y_i=0$ of $E$ over $U_i$'s.
Since the two ideals 
$\mathcal O_{U_i}x_i+\mathcal O_{U_i}y_i$ 
and $\mathcal O_{U_j}x_j+\mathcal O_{U_j}y_j$ define
the same ideal over $U_i \cap U_j$, 
there exist elements
$b_{ij}$ and $c_{ij}$ of $\mathcal O_{U_i \cap U_j}$ satisfying
$x_i = b_{ij} y_j + c_{ij} x_j$.
Then we have
$$
f_{ij}
= x_i s_j' - x_j s_i'
= (x_i - c_{ij} x_j) s_j' + (c_{ij} s_j' - s_i')x_j
= (b_{ij}y_j) s_j' + (c_{ij} s_j' - s_i')x_j
$$
and 
$$
\bar f_{ij}=\ov{(b_{ij}y_j)s_j'}+\ov{(c_{ij}s_j'-s_i')x_j}
\qquad
\mbox{in}
\qquad
\Gamma(S \cap U_i \cap U_j, \mathcal O_S)
$$
since $(b_{ij}y_j)s_j'$ belongs to $\mathcal O_{U_i\cap U_j}$.
Hence we have 
$$
d_S(v)_{ij}
= \dfrac{\bar f_{ij}}{\bar x_i \bar x_j}
= \dfrac{\ov{b_{ij}y_j}}{\bar x_i} \cdot \dfrac{\ov{s_j'}}{\bar x_j}
+ \dfrac{\ov{c_{ij}s_j'-s_i'}}{\bar x_i}
= \dfrac{\ov{b_{ij}y_j}}{\bar x_i} \cdot v
+ \dfrac{\ov{c_{ij}s_j'-s_i'}}{\bar x_i}
$$
in $\Gamma(S\cap U_i \cap U_j, \mathcal O_S(2E))$,
where $\ov{b_{ij}y_j}$ is 
a section of $\mathcal O_S(-S)\simeq {N_{S/V}}^{\vee}$ 
over $S \cap U_i \cap U_j$.
Since
$$
\ov{c_{ij}s_j'} - \ov{s_i'} 
= \bar c_{ij} \bar x_j v - \bar x_i v
=0
\qquad
\mbox{in}
\qquad
\Gamma(S \cap U_i \cap U_j, \mathcal O_S(S)),
$$
$c_{ij}s_j' - s_i' \in \mathcal O_{U_i \cap U_j}(S)$ is contained
in $\mathcal O_{U_i \cap U_j}$.
Hence $\ov{c_{ij}s_j'-s_i'}/\bar x_i$ is contained in 
$\Gamma(S\cap U_i \cap U_j,\mathcal O_S(E))$.
On the other hand, the restriction of the $1$-cochain 
$\left\{ \ov{b_{ij}y_j}/\bar x_i \right\}_{i,j \in I}$ to $E$ is a cocycle and 
represents the extension class $\mathbf e \in H^1(\mathcal O_E(-S+E))$ 
of the exact sequence \eqref{ses:normal bundle of line}.
Therefore $d_S(v)\big{\vert}_E$ is equal to
$\mathbf e \cup (v\big{\vert}_E)=\partial(v\big{\vert}_E)$,
which shows (2).
\QED

\subsection{Computation of obstructions}\label{computation}

The purpose of this subsection is the proof of Theorem~\ref{thm:main thm1}.
Let $v$ be a global section of $H^0(N_{S/V}(E))$ which satisfies 
the inequality \eqref{ineq:intersection}.
Let $\mathbf k_E=\mathbf k_{E,S}$ and 
$\mathbf k_C=\mathbf k_{C,S}$ be the extension classes of 
the exact sequences
$$
0 \longrightarrow \mathcal O_S(-E) 
\longrightarrow \mathcal O_S
\longrightarrow \mathcal O_E
\longrightarrow 0
\quad {\rm and} \quad
0 \longrightarrow \mathcal O_S(-C) 
\longrightarrow \mathcal O_S
\longrightarrow \mathcal O_C
\longrightarrow 0
$$
on $S$, respectively. 
We regard $v \big{\vert}_C$ (resp. $v \big{\vert}_E$) 
as a global section of $N_{S/V}\big\vert_C$ 
(resp. $N_{S/V}(E-C)\big\vert_E$).
Then we have the following:
\begin{lem}\label{lem:new}
$$
v \big{\vert}_C \cup \mathbf k_C = v \big{\vert}_E \cup \mathbf k_E
\qquad
\mbox{in}
\qquad
H^1(N_{S/V}(-C)).
$$
\end{lem}

\Proof
We have the following commutative diagram of $\mathcal O_S$-modules
$$
\left[\begin{array}{ccccccccc}
0 & \rightarrow & \mathcal O_S(-C-E) & \rightarrow & \mathcal O_S(-C)\oplus\mathcal O_S(-E) &\rightarrow& \mathcal I_{C\cap E} &\rightarrow & 0\\
 & & \mapdown{} && \mapdown{} & & \mapdown{|_C} && \\
0 & \rightarrow & \mathcal O_S(-C-E) & \rightarrow & \mathcal O_S(-E) &\rightarrow& \mathcal O_C(-E) &\rightarrow & 0
\end{array}\right] \otimes N_{S/V}(E)
$$
whose first row is the Koszul complex of $C \cap E$. 
By \eqref{ineq:intersection} the global section $v$ belongs to $H^0(\mathcal I_{C\cap E} \otimes N_{S/V}(E))$.
By the commutativity, the coboundary map 
$$
H^0(\mathcal I_{C\cap E}\otimes N_{S/V}(E))
\longrightarrow H^1(N_{S/V}(-C))
$$
of the first row is equal to $(\cup \mathbf k_C)\circ \big{\vert}_C $ and 
similarly to $ (\cup \mathbf k_E) \circ \big{\vert}_E$.
\QED

We need to consider the relation between the two maps $d_C$ and $d_S$ 
allowing pole along $E$. The diagram \eqref{diag:d-map} becomes 
the {\it partially commutative} diagram

\begin{equation}\label{diag:d-map with pole}
\begin{array}{ccccc}
v & \in & H^0(N_{S/V}(E)) & \mapright{d_S} & 
H^1(\mathcal O_S(2E)) \\
&& & & \mapdown{|_C}\\
&& \mapdown{|_C} && H^1(\mathcal O_C(2Z)) \\
&& && \mapdown{H^1(\iota)} \\
&& H^0(N_{S/V}(E)\big{\vert}_C)& & 
 H^1({N_{C/V}}^{\vee}\otimes N_{S/V}(2Z))\\
&& \cup &&\uparrow\\
u & \in & H^0(N_{S/V}\big{\vert}_C) & \mapright{d_C} & 
 H^1({N_{C/V}}^{\vee}\otimes N_{S/V}),
\end{array}
\end{equation}
where $Z$ is the scheme-theoretic intersection $C\cap E$.
In other words, the commutativity holds only for 
$u \in H^0(N_{S/V}\big{\vert}_C)$ which has a lift $v \in H^0(N_{S/V}(E))$.
More precisely, for such a pair $u$ and $v$, we have
\begin{equation}\label{eq:d-map with pole}
\overline{d_C (u)} = H^1(\iota)(d_S(v)\big{\vert}_C).
\end{equation}
Here $\overline{\ast}$ denotes the image of $\ast \in H^1(C, \mathcal F)$ 
(resp. $\ast \in H^1(S, \mathcal F)$) in $H^1(C, \mathcal F(2Z))$ 
(resp. $H^1(S, \mathcal F(2E))$), where $\mathcal F$ 
is a vector bundle on $C$ (resp. $S$).

\paragraph{Proof of Theorem~\ref{thm:main thm1}}
Let $\alpha \in H^0(N_{C/V})$ be as in the theorem.
It suffices to show that the exterior component $\ob_S(\alpha)$ 
is nonzero in $H^1(N_{S/V}\big{\vert}_C)$.
In fact, we show that its image $\overline{\ob_S(\alpha)}$ 
in $H^1(N_{S/V}(2E)\big{\vert}_C)$ is nonzero.
The following generalizes Lemma~\ref{lem:exterior component2} 
under the circumstances:

\medskip\noindent
{\bf Step 1} \quad 
$$
\overline{\ob_S(\alpha)}=d_{S}(v)\big{\vert}_C \cup \pi_S(\alpha).
$$

\smallskip
\Proof
By Lemma~\ref{lem:exterior component}, 
$\overline{\ob_S(\alpha)}$ is equal to the cup product
$\overline{d_C(\pi_S(\alpha))} \cup \alpha$.
By \eqref{eq:d-map with pole} $\overline{d_C (\pi_S(\alpha))}$ 
is equal to $H^1(\iota)(d_S(v)\big{\vert}_C)$.
The rest of the proof is same as that of 
Lemma~\ref{lem:exterior component2}. By the commutative diagram
$$
\begin{array}{ccccc}
H^1({N_{C/V}}^{\vee}\otimes N_{S/V}(2E)) & \times 
 & H^0(N_{C/V}) & \overset{\cup}{\longrightarrow} & H^1(N_{S/V}(2E)\big{\vert}_C)\\
\mapup{H^1(\iota)} && \mapdown{\pi_S} && \Vert \\
H^1(\mathcal O_C(2Z)) & \times & H^0(N_{S/V}\big{\vert}_C) 
 & \overset{\cup}{\longrightarrow} & H^1(N_{S/V}(2E)\big{\vert}_C),
\end{array}
$$
we have the required equation.
\QED

We relate $\ob_S(\alpha)$ with a cohomology class on $E$ by
Lemma~\ref{lem:new}:

\medskip\noindent
{\bf Step 2} \quad 
$$
\overline{\ob_S(\alpha)} \cup \mathbf k_C 
=(d_S(v)\big{\vert}_E \cup v\big{\vert}_E) \cup \mathbf k_E 
$$
in $H^2(N_{S/V}(2E-C))$. 

\smallskip
\Proof
Since
$$
\begin{array}{ccccc}
H^1(\mathcal O_C(2Z)) & \times & H^0(N_{S/V}\big{\vert}_C) 
 & \overset{\cup}{\longrightarrow} & H^1(N_{S/V}(2E)\big{\vert}_C)\\
\mapup{|_C} && \mapid && \mapid \\
H^1(\mathcal O_S(2E)) & \times & H^0(N_{S/V}\big{\vert}_C) 
 & \overset{\cup}{\longrightarrow} & H^1(N_{S/V}(2E)\big{\vert}_C)\\
\mapid && \mapdown{\cup \, \mathbf k_C} && 
\mapdown{\cup \, \mathbf k_C} \\
H^1(\mathcal O_S(2E)) & \times & H^1(N_{S/V}(-C))
 & \overset{\cup}{\longrightarrow} & H^2(N_{S/V}(2E-C))
\end{array}
$$
is commutative, we have a commutative diagram
\begin{equation}\label{diag:res-coboundary1}
\begin{array}{ccccc}
 && \pi_S(\alpha) && \overline{\ob_S(\alpha)}\\ [-10pt]
&& \rotatebox{-90}{$\in$} && \rotatebox{-90}{$\in$} \\ [4pt]
H^1(\mathcal O_C(2Z)) & \times & H^0(N_{S/V}\big{\vert}_C) 
 & \overset{\cup}{\longrightarrow} & H^1(N_{S/V}(2E)\big{\vert}_C)\\
\mapup{|_C} && \mapdown{\cup \, \mathbf k_C} && 
\mapdown{\cup \, \mathbf k_C} \\
H^1(\mathcal O_S(2E)) & \times & H^1(N_{S/V}(-C))
 & \overset{\cup}{\longrightarrow} & H^2(N_{S/V}(2E-C)). \\
\rotatebox{90}{$\in$} && && \\ [-4pt] 
d_S(v) && && \\ 
\end{array}
\end{equation}
Hence we have 
$$
\overline{\ob_S(\alpha)} \cup \mathbf k_C
=(d_{S}(v)\big{\vert}_C \cup \pi_S(\alpha)) \cup \mathbf k_C
= d_{S}(v) \cup (\pi_S(\alpha) \cup \mathbf k_C)
$$
in $H^2(N_{S/V}(2E-C))$ by Step 1.
There exists a commutative diagram
\begin{equation}\label{diag:res-coboundary2}
\begin{array}{ccccc}
&& v\big{\vert}_E && \\ [-10pt]
&& \rotatebox{-90}{$\in$} && \\ [4pt]
H^1(\mathcal O_E(2E)) & \times & H^0(N_{S/V}(E-C)\big{\vert}_E)
 & \overset{\cup}{\longrightarrow} & H^1(N_{S/V}(3E-C)\big{\vert}_E) \\
\mapup{|_E} && \mapdown{\cup \, \mathbf k_E} && 
\mapdown{\cup \, \mathbf k_E} \\
H^1(\mathcal O_S(2E)) & \times & H^1(N_{S/V}(-C))
 & \overset{\cup}{\longrightarrow} & H^2(N_{S/V}(2E-C)), \\
\rotatebox{90}{$\in$} && && \\ [-4pt] 
d_S(v) && && \\ 
\end{array}
 \end{equation}
similar to \eqref{diag:res-coboundary1}.
Therefore by Lemma~\ref{lem:new}, we have
$$
d_{S}(v) \cup (\pi_S(\alpha) \cup \mathbf k_C)
=d_{S}(v) \cup (v\big{\vert}_E \cup \mathbf k_E)
=(d_{S}(v)\big{\vert}_E \cup v\big{\vert}_E) \cup \mathbf k_E.
$$
Thus we obtain the equation required.
\QED

\medskip\noindent
{\bf Step 3} \quad 
Since $d_{S}(v)\big{\vert}_E=\partial (v\big{\vert}_E)$
by Proposition~\ref{prop:key proposition} (2),
we obtain $d_{S}(v)\big{\vert}_E \ne 0$ by the assumption (b).
Since $N_{S/V}(E-C)\big{\vert}_E$ is trivial by the assumption (a),
we have $d_{S}(v)\big{\vert}_E \cup v\big{\vert}_E \ne 0$
in $H^1(N_{S/V}(3E-C)\big{\vert}_E) \simeq H^1(\mathcal O_E(2E))$.
Consider the coboundary map
$$
\cup\, \mathbf k_E: H^1(N_{S/V}(3E-C)\big{\vert}_E) \longrightarrow
H^2(N_{S/V}(2E-C)),
$$
which appears in \eqref{diag:res-coboundary2}.
By the Serre duality, it is the dual of the restriction map
$$
H^0(S,C+K_V\big\vert_S -2E) \overset{|_E}{\longrightarrow} H^0(E,(C+K_V-2E)\big\vert_E),
$$
which is surjective by the assumption (c).
Hence the coboundary map $\cup\, \mathbf k_E$ is injective.
Therefore we obtain
$d_{S}(v)\big{\vert}_E \cup v\big{\vert}_E \cup \mathbf k_E\ne 0$
and hence by Step 2 we conclude that
$\overline{\ob_S(\alpha)}\ne 0$ in $H^1(N_{S/V}(2E)\big{\vert}_C)$.

\medskip

Thus we have proved Theorem~\ref{thm:main thm1}.

\section{Application to Hilbert schemes}\label{hilbertscheme}

In this section, we apply the result of the previous section 
to prove Theorem~\ref{thm:main thm2}.
We generalize Mumford's example (Example~\ref{ex:pathology}) and show that
for many uniruled $3$-folds $V$, their Hilbert schemes 
$\Hilb^{sc} V$ contain similar generically non-reduced components.

\subsection{Dichotomy}\label{dichotomy}

We explain the detail of Example~\ref{ex:cubic 3-fold}, which 
is a prototype of the non-reduced components constructed in \S\ref{proof}.
It is simpler than Mumford's example in applying Theorem~\ref{thm:main thm1}.
Let $C,E,S,V,h$ and $W$ be as in Example~\ref{ex:cubic 3-fold}.
Then by $(C. E)=(2h +2E. E)=0$, the intersection $C \cap E$ is empty,
which is the main reason for the simplicity.
By adjunction, we have 
$\mathcal O_S(-K_S) \simeq \mathcal O_S(h)\simeq N_{S/V}$.
By adjunction again, $N_{S/V}\big{\vert}_C$ and $N_{C/S}$ are isomorphic to
$\mathcal O_C(K_C)$ and $\mathcal O_C(2K_C)$, respectively.
By the exact sequence \eqref{ses:normal bundle of curve},
$h^0(N_{C/V})$ is equal to 
$$
h^0(N_{C/S})+h^0(N_{S/V}\big{\vert}_C)
=h^0(2K_C)+h^0(K_C)=12+5=17.
$$
Hence the tangent space $t_{W,C}$ of $W$ at $[C]$
is of codimension one in the tangent space $H^0(N_{C/V})$ 
of $\Hilb^{sc} V$.
We have only the following two possibilities ({\it i.e}. dichotomy):
\begin{enumerate}
 \renewcommand{\labelenumi}{{\rm [\Alph{enumi}]}}
 \item $W$ is an irreducible component of $(\Hilb^{sc} V)_{\red}$.
  Moreover $\Hilb^{sc} V$ is generically non-reduced along $W$.
 \item There exists an irreducible component $W'$ of $\Hilb^{sc} V$
  which contains $W$ as a proper closed subset.
  $\Hilb^{sc} V$ is generically smooth along $W$.
\end{enumerate}
We prove that the case [B] does not occur.

\begin{prop}\label{prop:cubic 3-fold}
 The Hilbert scheme of smooth curves on a smooth cubic $3$-fold $V$
 contains a generically non-reduced component $\tilde W$
 of dimension $16$ such that $(\tilde W)_{\red}=W$.
\end{prop}

\Proof 
Consider the exact sequence
$$
0 \longrightarrow N_{S/V}
\longrightarrow N_{S/V}(E)
\longrightarrow N_{S/V}(E)\big{\vert}_E
\longrightarrow 0.
$$
Since $H^1(N_{S/V})\simeq H^1(\mathcal O_S(h))=0$, 
there exists a rational section $v$ of $N_{S/V}$
having a pole along $E$.
By $H^1(N_{C/S})\simeq H^1(2K_C)=0$ and \eqref{ses:normal bundle of curve}, 
there exists a first order infinitesimal deformation 
$\tilde C \subset V\times \Spec k[t]/(t^2)$ of $C \subset V$
whose exterior component $\pi_S(\alpha)$ coincides with
$v\big{\vert}_C$.
Since $S$ is general, so is $E$.
Hence, by \cite[Proposition 1.3]{Iskovskih77}, the normal bundle 
$N_{E/V}$ is trivial.
Moreover it is easily checked that 
the other conditions of Theorem~\ref{thm:main thm1} are satisfied.
For example, since $v\big{\vert}_E \ne 0$ and
the coboundary map $\partial$ is injective, we have (a). 
Since $C \sim -K_V\big{\vert}_S+2E$, we have $\Delta=0$ by definition.
Therefore we have (b).
Thus we conclude that $\tilde C$ is obstructed by the theorem.
This implies [A].
\QED

\begin{rmk}
 By a similar method, we can show the non-reducedness of $\Hilb^{sc}
 \mathbb P^3$ along $W=W^{56}$
 for Mumford's example in arbitrary characteristic.
\end{rmk}

\subsection{$S$-maximal family of curves}\label{maximal}

Let $V$ be a smooth projective $3$-fold 
and let $S$ be a smooth surface in $V$.
We introduce the notion of $S$-maximal families, which is 
analogous to $s$-maximal irreducible subsets defined in 
\cite{Kleppe} for $\Hilb^{sc} \mathbb P^3$.
We assume that the Hilbert scheme $\Hilb V$ of $V$ is nonsingular at $[S]$.
Let $U_S$ be the irreducible component passing through $[S]$ and
let 
$$
V \times U_S \supset \mathcal S \overset{p_2}\longrightarrow U_S
$$
be the universal family over $U_S$.
Let $C$ be a smooth curve on $S$ and
assume that $\Hilb^{sc} S$ is nonsingular 
of expected dimension (=$\chi(N_{C/S})$) at $[C]$ 
({\it i.e.} $H^1(N_{C/S})=0$).
Then the Hilbert scheme $\Hilb^{sc} \mathcal S$,
which is same as 
the relative Hilbert scheme of $\mathcal S/U_S$ is nonsingular at $[C]$.
The first projection $p_1 : \mathcal S \rightarrow V$ induces
the morphism $\Hilb^{sc} \mathcal S \rightarrow \Hilb^{sc} V$.
Let $\mathcal W_{S,C}$ be the irreducible component of $\Hilb^{sc} \mathcal S$
passing through $[C]$.
We call the image of $\mathcal W_{S,C}$ in $\Hilb^{sc} V$
the {\em $S$-maximal family of curves} containing $C$ 
and denote it by $W_{S,C}$.
We illustrate $W_{S,C}$ by the diagram
$$
\begin{array}{cccccc}
\mathcal W_{S,C} &\subset & \Hilb^{sc} \mathcal S & & \mathcal S 
& ({\rm universal \ family}) \\
\downarrow && \downarrow & \searrow & \downarrow & \\
W_{S,C} & \subset & \Hilb V & \supset & U_S. & 
\end{array}
$$
There exists a commutative diagram
$$
\begin{array}{ccccccc}
 0 \longrightarrow 
 & N_{C/S} & 
 \longrightarrow & N_{C/V} &
 \overset{\pi_S}{\longrightarrow}
 & N_{S/V}\big{\vert}_C & 
 \longrightarrow 0 \\
 & || && \uparrow && \uparrow&\\
 0 \longrightarrow 
 & N_{C/S}&
 \longrightarrow & N_{C/\mathcal S} &
 \longrightarrow
 & \underbrace{N_{S/\mathcal S}\big{\vert}_C}_{\cong \, H^0(N_{S/V}) \otimes \mathcal O_C}& 
 \longrightarrow 0.
 \end{array}
$$
Here the two horizontal sequences of normal bundles are exact.
By the diagram, we have the following lemma.
\begin{lem}\label{lem:tangential map}
The cokernel (resp. kernel) of the tangential map
\begin{equation}\label{map:tangential map}
\kappa_{[C]}: H^0(N_{C/\mathcal S}) \longrightarrow H^0(N_{C/V}) 
\end{equation}
of $\mathcal W_{S,C} \rightarrow \Hilb^{sc} V$ at $[C]$ 
is isomorphic to that of the restriction map
$H^0(N_{S/V}) \rightarrow H^0(N_{S/V}\big{\vert}_C)$.
\end{lem}

\subsection{Construction of obstructed curves}\label{construction}

From now on we assume that the geometric genus 
$p_g(S)$ is zero and $H^1(N_{S/V})=0$.
Let $E$ be a $(-1)$-$\mathbb P^1$ on $S$ whose normal bundle 
$N_{E/V}$ is generated by global sections.
We denote by $\varepsilon: S \rightarrow F$ the blow-down of $E$ from $S$.

\begin{prop}\label{prop:delta}
Let $\Delta_1$ be a very ample divisor on $F$.
Then for each sufficiently large integer $n$, $\Delta_n := n\Delta_1$
 satisfies the following conditions: 
 \begin{enumerate}
\renewcommand{\labelenumi}{{\rm [\roman{enumi}]}}
 \setlength{\itemsep}{-3pt} 
 \item the linear system 
  $\Lambda_n:=|\varepsilon^* \Delta_n -K_V\big\vert_S +2E|$ on $S$
  contains a smooth connected member $C$,
 \item the restriction map  
	$\Lambda_n \cdots \rightarrow \Lambda_n\big{\vert}_E$ is surjective,
 \item $H^i(S,\varepsilon^* \Delta_n+E)=0$ for $i =1,2$,
 \item $H^1(S,\varepsilon^* \Delta_n-E)=0$ and
 \item $H^1(C,N_{C/S})=0$.
\end{enumerate}
\end{prop}
In the next subsection, we will show that $\Hilb^{sc} V$ 
is non-reduced in a neighborhood of the corresponding point $[C]$.

\smallskip
\Proof
We have $\chi(N_{E/V})=\deg (-K_V\big\vert_E) =\deg N_{E/V}+2$.
Since $N_{E/V}$ is generated by global sections,
we have $\deg N_{E/V}\ge 0$ and hence $\deg (-K_V\big\vert_E)\ge 2$.
We have $(D. E) \ge -1$ for 
$D=E$, $-E$, $-K_V\big\vert_S +E$ and $-K_V\big\vert_S +2E$.
By the lemma below, there exists an integer $m_1$ such that
for each $n \ge m_1$, all of the cohomology groups
$H^i(\varepsilon^* (n\Delta_1)+E)$ ($i=1,2$),
$H^1(\varepsilon^* (n\Delta_1)-E)$,
$H^1(\varepsilon^* (n\Delta_1)-K_V\big\vert_S +E)$ 
and $H^1(\varepsilon^* (n\Delta_1)-K_V\big\vert_S +2E)$ vanish.

Put $e:=\deg (-K_V\big{\vert}_S + 2E)\big{\vert}_E$. Then $e \ge 0$.
Suppose that $e=0$. Then there exists an integer $m_2$ such that
for each $n \ge m_2$, $n\Delta_1+\varepsilon_*(-K_V\big\vert_S +2E)$
is very ample on $F$. Hence by the Bertini theorem
({\it cf}.~\cite[Chap.~II, Theorem 8.18]{Hartshorne}),
the linear system $|n\Delta_1+\varepsilon_*(-K_V\big\vert_S +2E)|$ 
contains a smooth connected member.
Suppose that $e>0$. Then there exists an integer $m_2$
such that for each $n \ge m_2$,
$\Lambda_n=\left| \varepsilon^* (n\Delta_1)-K_V\big\vert_S +2E \right|$ 
is base point free and $\varepsilon^*(n\Delta_1)-K_V\big\vert_S +2E$ 
is ample. Then by the Bertini theorem again and 
\cite[Chap.~III Corollary 7.9]{Hartshorne},
$\Lambda_n$  contains a smooth connected member.

Assume that $n \ge \max\{m_1, m_2\}$  and let  $C$  be as in [i].
Then we obtain [ii], [iii] and [iv] by the choice of $m_1$.
Finally we prove [v]. Consider the exact sequence
$$
0 \longrightarrow \mathcal O_S
\longrightarrow \mathcal O_S(C)
\longrightarrow N_{C/S}
\longrightarrow 0,
$$
which induces
$$
H^1(S,\mathcal O_S(C)) \longrightarrow H^1(C,N_{C/S})
\longrightarrow H^2(S,\mathcal O_S).
$$
It follows from the choice of $m_1$ that
$H^1(C)\simeq H^1(\varepsilon^* (n\Delta_1)-K_V\big\vert_S +2E)=0$.
Since $p_g(S)=0$, we have $H^1(N_{C/S})=0$.
\QED

\begin{lem}\label{lem:vanish}
 Let $D$ be a divisor on $S$ with $(D. E) \ge -1$.
 Then there exists an integer $m_0$ such that
 for each $i >0$ and each $n \ge m_0$,
 we have $H^i(S,\varepsilon^* (n\Delta_1)+D)=0$.
\end{lem}
\Proof
By assumption, $D$ is linearly equivalent to 
$\varepsilon^*D'-jE$ for some divisor $D'$ on $F$, 
where $j=(D.E)\ge -1$. If $j \ge 0$, by the Serre vanishing theorem, 
there exists an integer $m_0$ such that
$H^i(S,\varepsilon^* (n\Delta_1)+D)
\simeq H^i(F,\mathfrak m_p^j(D'+n\Delta_1))=0$ for each $n \ge m_0$,
where $\mathfrak m_p$ is the maximal ideal at $p=\varepsilon(E)$.
If $j=(D.E)=-1$, then we have $H^i(\mathcal O_E(D\big{\vert}_E))=0$.
Hence by the exact sequence
$$
[0 \longrightarrow \mathcal O_S(D-E)
\longrightarrow \mathcal O_S(D)
\longrightarrow \mathcal O_E(D\big{\vert}_E)
\longrightarrow 0] \otimes \mathcal O_S(\varepsilon^* (n\Delta_1)),
$$
the assertion follows from the case $j=0$.
\QED

\subsection{Proof of Theorem~\ref{thm:main thm2}}\label{proof}

Let $C \in \Lambda_n$ be as in [i] of Proposition~\ref{prop:delta} for a sufficiently large integer $n$. 
By assumption and [v],
$\Hilb V$, $\Hilb^{sc} S$ and $\Hilb^{sc} \mathcal S$ are 
nonsingular (of expected dimension) at $[S]$, $[C]$ and $[(C,S)]$,
respectively.
We consider the $S$-maximal family $W_{S,C}$ of curves containing $C$
({\it cf}. \S\ref{maximal}).
Let $\kappa_{[C]}: H^0(N_{C/\mathcal S}) \rightarrow H^0(N_{C/V})$
be the tangential map
of the morphism $\Hilb^{sc} \mathcal S \rightarrow \Hilb^{sc} V$ 
({\it cf}. \eqref{map:tangential map}).
Then we have the following:

\medskip\noindent
{\bf Claim 1} \quad 
$\dim\coker\kappa_{[C]}=1$.
\medskip

\Proof
By Lemma~\ref{lem:tangential map},
the cokernel is isomorphic to $H^1(N_{S/V}(-C))$ since $H^1(N_{S/V})=0$.
Since 
$N_{S/V}(-C) \sim K_S-K_V\big{\vert}_S -C \sim K_S-\varepsilon^* \Delta_n-2E$,
we have $H^1(N_{S/V}(-C)) \simeq H^1(\varepsilon^* \Delta_n+2E)^{\vee}$ 
by the Serre duality.
Let us consider the exact sequence
$$
0 \longrightarrow \mathcal O_S(\varepsilon^* \Delta_n+E)
\longrightarrow \mathcal O_S(\varepsilon^* \Delta_n+2E)
\longrightarrow \mathcal O_E(2E)
\longrightarrow 0.
$$
By the condition [iii], we have 
$H^1(\varepsilon^* \Delta_n+2E)\simeq H^1(\mathcal O_E(2E))\simeq k$.
\QED

If $\alpha \in \im \kappa_{[C]}$, then 
the corresponding first order infinitesimal deformation $\tilde C$
of $C \subset V$ is realized as a member of $W_{S,C}$.
Hence by Claim 1, the same dichotomy between [A] and [B]
in \S\ref{dichotomy} holds for $W:=W_{S,C}$.

\medskip\noindent
{\bf Claim 2} \quad 
If $\alpha \not\in \im \kappa_{[C]}$, then $\tilde C$ is obstructed.
\medskip

\Proof
We show that the exterior component 
$\pi_S(\alpha) \in H^0(N_{S/V}\big{\vert}_{C})$ 
of $\alpha$ lifts to a rational section 
$v$ of $N_{S/V}$ having a pole of order one along $E$.
Consider a commutative diagram
$$
\begin{array}{ccccc}
H^0(N_{S/V}) & \overset{|_{C}}{\longrightarrow} & 
H^0(N_{S/V}\big{\vert}_{C}) & \twoheadrightarrow & 
\phantom{\simeq k} H^1(N_{S/V}(-C)) \simeq k \\
\bigcap && \bigcap && \downarrow \\
H^0(N_{S/V}(E)) & \longrightarrow & H^0(N_{S/V}(E)\big{\vert}_{C})
& \longrightarrow & H^1(N_{S/V}(E-C)). \\
\end{array}
$$
By the Serre duality and the condition [iii],
we have $H^1(N_{S/V}(E-C))\simeq H^1(\varepsilon^* \Delta_n+E)^{\vee}=0$.
Hence there exists $v \in H^0(N_{S/V}(E))$
such that $v\big{\vert}_{C}=\pi_S(\alpha)$.
By the choice of $\alpha$, $v$ is not contained in $H^0(N_{S/V})$.

Now we check that the two assumptions (a) and (b) of 
Theorem~\ref{thm:main thm1} are satisfied.
First we consider (a). Since $v\big{\vert}_{C}=\pi_S(\alpha)$ 
is contained in $H^0(N_{S/V}\big{\vert}_{C})$,
we have $(v)_0 \cap E \ge C \cap E$ as divisor on 
$E \simeq \mathbb P^1$. Note that
$$
m:=(C. E)
=(-K_V\big{\vert}_S+2E. E)
=\deg (-K_V\big\vert_E)-2
=\deg N_{E/V}=\deg N_{S/V}(E)\big{\vert}_E.
$$
By the degree reason, 
we have $(v)_0 \cap E=C \cap E$.
Since $C$ is a general member of $\Lambda_n$, by the condition [ii],
$C$ meets $E$ at general $m$ points on $E$.
Hence $v\big{\vert}_E$ is a general global section of
$N_{S/V}(E)\big{\vert}_E$.
Therefore we have (a) by Lemma~\ref{lem:nonzero} below.
Consider the exact sequence
$0 \rightarrow \mathcal O_S(\varepsilon^* \Delta_n -E) 
\rightarrow \mathcal O_S(\varepsilon^* \Delta_n) 
\rightarrow \mathcal O_E
\rightarrow 0$ for (b). It follows from [iv] that the restriction map
$H^0(S,\varepsilon^* \Delta_n)\rightarrow H^0(E,\mathcal O_E)$ is surjective.
By Theorem~\ref{thm:main thm1}, $\tilde C$ is obstructed.
\QED

\begin{lem}\label{lem:nonzero}
 Let $\partial$ be the coboundary map in Theorem~\ref{thm:main thm1}.
 If $N_{E/V}$ is generated by global sections and
 $t$ is a general global section of $N_{S/V}(E)\big{\vert}_E$, 
 then the image $\partial (t)$ is nonzero in $H^1(\mathcal O_E(2E))$.
\end{lem}
\Proof
By assumption, we have $H^1(N_{E/V}(E))=0$ and hence
$\partial$ is surjective.
Since $t \in H^0(N_{S/V}(E)\big{\vert}_E)$ is general,
it is not contained in the kernel of the coboundary map
and hence $\partial(t)\ne 0$.
\QED

Therefore, by Claim 2, we have [A]
and hence we have proved Theorem~\ref{thm:main thm2}.

\subsection{Examples}\label{examples}

If $V$ is separably uniruled and $E\simeq \mathbb P^1$ sweeps out $V$,
then $N_{E/V}$ is generated by global sections.
\begin{ex}[$\car k=0$]\label{ex:uniruled}
The following $V, S$ and $E$ satisfy the assumption of 
Theorem~\ref{thm:main thm2}.
\begin{enumerate}
\item $V$ is a Fano $3$-fold whose anti-canonical class $-K_V$ is 
      a sum $H+H'$ of two ample divisors $H$ and $H'$.
      Then $|H|$  contains a smooth member  $S$.
      If $(V, H) \not \simeq 
      (\mathbb P^3, \mathcal O_{\mathbb P}(1)),
      (\mathbb P^3, \mathcal O_{\mathbb P}(2)),
      (Q^3, \mathcal O_Q(1))$, 
      then $S$ contains  a $(-1)$-$\mathbb P^1$  $E$.
  
\item $V$ has a $\mathbb P^1$-bundle structure  
      $\pi : V \rightarrow F$ 
      over a smooth surface $F$ with $p_g=0$  in Zariski topology.
      Let $S_1 \subset V$ be a rational section of  $\pi$ 
      and $S$ a smooth member of $|S_1+ \pi^*A|$ 
      for a sufficiently ample divisor $A$ on $F$.
      Then $\pi\vert_S : S \rightarrow F$  is birational but not isomorphic.
      Hence $S$  contains  a $(-1)$-$\mathbb P^1$  $E$.

\item $V$ has a $\mathbb P^2$-bundle structure $\pi : V \rightarrow C$ over a smooth curve $C$. 
Let  $\mathcal O(1)$  be a relative tautological line bundle 
and $S$ a smooth member of  $|\mathcal O(2)(\pi^*A)|$ for a sufficiently ample divisor  $A$ on $C$.
Then $\pi\vert_S : S \to C$  is  a conic bundle with a reducible fiber.
An irreducible component $E$ of a reducible fiber is a $(-1)$-$\mathbb P^1$.
\end{enumerate}
\end{ex}

\begin{rmk}
 Deformations of the curve $C \subset \mathbb P^4$ 
 in Example~\ref{ex:cubic 3-fold} in a smooth quintic $3$-fold 
 $V_5 \subset \mathbb P^4$ 
 was studied as Voisin's example in Clemens-Kley~\cite{CK}.
 They showed that the Hilbert scheme $\Hilb^{sc}_{8,5} V_5$ 
 has an embedded component at $[C]$.
\end{rmk}

\begin{rmk}\label{rmk:murphy}
 Recently Vakil~\cite{V} has shown that various moduli spaces, 
 including the Hilbert schemes of space curves 
 (but not including $\Hilb^{sc} \mathbb P^3$), satisfy ``Murphy's law'',
 {\it i.e.}, every singularity type of finite type over $\mathbb Z$
 appears on the moduli spaces.
\end{rmk}

\section{Application to Hom schemes}\label{hom scheme}

Let $V$ be a smooth projective variety and $X$ a (smooth projective) curve.
It is well known that the Zariski tangent space of $\Hom(X,V)$
at $[f]$ is isomorphic to $H^0(X,f^* T_V)$ 
and the following dimension estimate holds:
\begin{equation}\label{neq:bound}
 \chi(f^* T_V) \le \dim_{[f]} \Hom(X,V) \le \dim H^0(X,f^* T_V),
\end{equation}
where $T_V$ is the tangent bundle of $V$.
The lower bound $\chi(f^* T_V)=\deg f^*(-K_V) +n(1-g)$ is called the 
{\it expected dimension}, where $n = \dim V$ and $g$ is the genus of $X$.

In this section we prove Theorem~\ref{thm:main thm}.
First we recall a simple unprojection.

\begin{lem}\label{lem:unprojection}
A smooth cubic surface $S \subset \mathbb P^3$ is isomorphic to the image of a smooth quartic del Pezzo surface $F= (2) \cap (2) \subset \mathbb{P}^4$ by a projection from a point on $F$.
\end{lem}
\Proof
As is well known $S$ contains a line $E$.
Choose a homogeneous coordinates $(x_1: x_2: x_3: x_4)$ 
of $\mathbb{P}^3$ such that $E$ is defined by $x_1=x_2=0$.
Then the equation of $S$ is $x_1q(x)+x_2q'(x)=0$ for quadratic forms 
$q(x)$ and $q'(x)$. $S$ is the image of the quartic del Pezzo surface 
$F:q(x)+x_2y=q'(x)-x_1y=0$ in the projective 4-space $\mathbb P^4$ 
with coordinate $(x: y)$ from the point $(0:0:0:0:1) \in F$.
\QED

Let $X$ be a general curve of genus 5.
The canonical model of $X$, that is, the image of 
$X \stackrel{K_X}\hookrightarrow \mathbb P^4$, is a 
general smooth complete intersection $q_1=q_2=q_3=0$ of three quadrics.
Let $p$ be a general point of the ambient space $\mathbb P^4$ and 
$F_p$ the intersection of two quadrics $q$ and $q'$ 
which belong to the net of quadrics $\langle q_1, q_2, q_3 \rangle_k$ 
and pass through $p$.
We denote the blow-up of $F_p$ at $p$ by $\pi_p: S_p \rightarrow F_p$.
Then we have a commutative diagram:
\begin{equation}\label{diag:projection}
\begin{array}{cccccl}
X & \subset & F_p & \subset & \mathbb P^4 & \\
|| && \smapup{\pi_p} & & \rmapdown{\Pi_p} & \\
C & \subset & S_p & \subset & \mathbb P^3, &
\end{array}
\end{equation}
where $C$ is the inverse image of $X$ in $S_p$ and $\Pi_p$ 
is the projection from $p \in F_p \setminus X$.
Since $X$ belongs to the linear system $|-2K_F|$ on $F_p$, 
$C$ belongs to $|\pi_p^*(-2K_F)| = |-2K_S+2E|$ on $S_p$, where $E$ 
is the exceptional curve of $S_p \rightarrow F_p$.
Since $q, q'$ and $p \in F_p$ are general, 
$S_p$ is a general cubic surface by Lemma~\ref{lem:unprojection}.

Let $\tilde W$ be the generically non-reduced component of 
$\Hilb^{sc} V_3$ in Proposition~\ref{prop:cubic 3-fold}.
Let $\varphi: \tilde W \rightarrow {\mathfrak M}_5$ be the 
classification morphism of $\tilde W$ to the moduli space 
of curves of genus 5.
The fiber $\varphi^{-1}([X])$ at the point $[X] \in {\mathfrak M}_5$ 
is isomorphic to an open subscheme of $\Hom(X,V_3)$.
We show that its Zariski closure $T$ in $\Hom(X,V_3)$ 
satisfies the requirement of Theorem~\ref{thm:main thm}.

In the Fermat case, every general cubic surface is isomorphic to 
a hyperplane section of $V_3^{\rm Fermat}$ by the following theorem,
for which we need $\car k=0$.

\smallskip
{\bf Sylvester's pentahedral theorem} (\cite[Chap.~\S 84]{Segre},~\cite{DvG})
{\it Every general cubic form $F(y_0, y_1, y_2, y_3)$ of $4$ variables is 
a sum $\sum_{i=0}^4 l_i(y_0, y_1, y_2, y_3)^3$ of the cubes of $5$ 
linear forms.}

\smallskip
Hence $S_p$ is isomorphic to a hyperplane section of $V_3^{\rm Fermat}$ 
and $C$ is a general member of $\tilde W^{\rm Fermat}$.
Therefore, the classification morphism 
$\varphi^{\rm Fermat}: \tilde W^{\rm Fermat} \rightarrow {\mathfrak M}_5$
is dominant by the diagram (\ref{diag:projection}) and
the fiber $T^{\rm Fermat}$ is of dimension $4$.
Since ${\mathfrak M}_5$ is generically smooth, $T^{\rm Fermat}$ 
is generically non-reduced.

Theorem~\ref{thm:main thm} for a general $V_3$ follows from the Fermat 
case by the upper semi-continuity theorem on fiber dimensions.
\QED

Theorem~\ref{thm:main thm} holds true for a smooth $V_3$ if the answer 
to following is affirmative:
\begin{prob}
Is the classification map
$$
(\mathbb P^4)^\vee \dashrightarrow \mathfrak M_{\rm cubic},
\qquad [H] \mapsto [H \cap V_3]
$$
dominant for a smooth cubic 3-fold $V_3 \subset \mathbb P^4$?
\end{prob}

\bigskip
\noindent  Research Institute for Mathematical Sciences \\
Kyoto University \\
Kyoto 606-8502\\
Japan \\
 {\it e-mail addresses} : mukai@kurims.kyoto-u.ac.jp and nasu@kurims.kyoto-u.ac.jp

\end{document}